-----------------------------------------------------------
\documentstyle[12pt]{article}
\title{{\normalsize{{\hskip 8.5cm} BIHEP-TH-95-25}}\\[-7mm]
{\normalsize{{\hskip 8.5cm} Sept.~~~~1995}}\\
{$L^{\pm}$ Operators and}\\
{Quantum Enveloping Algebras}
\thanks{This work was supported by the National Natural Science
        Foundation of China and Grant No.LWTZ-1298 of Chinese
        Academy of Sciences.}}
\author{{Mi Xie}\\
{\small Graduate School, Chinese Academy of Sciences}\\
        {Bai-Qi Jin$~~~~~~$Zhong-Qi Ma}\\
{\small Institute of High Energy Physics, Chinese Academy of Sciences}}
\date{}

\begin{document}

\maketitle

\begin{center}
\begin{minipage}{12cm}
\vspace{4cm}
\begin{center}
{\bf Abstract}
\end{center}

\vspace{2mm}
{The relations between $L^{\pm}$
operators and the generators in the quantum enveloping algebras
are studied. The $L^{\pm}$ operators for $U_{q}A_{N}$ and
$U_{q}G_{2}$ algebras are explicitly expressed by the generators
as examples.}

\end{minipage}
\end{center}

\newpage
\section*{I. Introduction}

\hskip 20pt
Recently, the theories of quantum enveloping algebras and quantum
groups have drawn wide concerns, and made great progress $^{[1,2]}$.
For a generic $q$, there have been standard methods for calculating
the highest weight representations and Clebsch-Gordan coefficients
of quantum enveloping algebras $^{[3,4]}$. Based on them, a series
of solutions $R_{q}$ of the simple Yang-Baxter equation (without
spectral parameter) and solutions $R_{q}(x)$ of the Yang-Baxter
equation $^{[5]}$ were found. In terms of the quantum double, the
general form of universal ${\cal R}$ matrix were obtained in principle
$^{[3]}$. For $q$ being a root of unity, the representation theory was
studied, too $^{[6]}$. The representation of quantum affine
algebras and their applications to the solvable lattice models
has been made great advance $^{[7]}$. On the other hand, following
the general idea of Connes on the non-commutative geometry $^{[8]}$,
Woronowicz $^{[9]}$ elaborated the framework of the non-commutative
differential calculus. He introduced the bimodule over the quantum
groups and presented varies theorems concerning the differential
forms and exterior derivative. The $q$-deformed gauge theory was
studied from several different viewpoints $^{[8-10]}$.

A quantum group is introduced as a Hopf algebra ${\cal A}=Fun_{q}(G)$
which is both non-commutative and non-cocommutative. It is
the continuous deformation of the Hopf algebra of the functions
on a Lie group, and is freely generated by the non-commutative
matrix elements $T^{a}_{~b}$:
$$(R_{q})^{ca}_{~rs}T^{r}_{~b}T^{s}_{~d}~=~
T^{a}_{~r}T^{c}_{~s}(R_{q})^{sr}_{~bd} \eqno (1.1) $$

\noindent
where the matrix $T$ takes value in the minimal representation, and
$R_{q}$ is the solution of the simple Yang-Baxter equation corresponding
to the minimal representation.

Faddeev-Reshetikhin-Takhtajan $^{[11]}$ introduced two linear functionals
$L^{\pm}$ which belong to the dual Hopf algebra ${\cal A'}$. They are
defined by their values on the elements $T^{a}_{~b}$:

$$\left(L^{+}\right)^{a}_{~b}\left(T^{c}_{~d}\right)~=~
\left(R_{q}^{-1}\right)^{ac}_{~bd},~~~~~
\left(L^{-}\right)^{a}_{~b}\left(T^{c}_{~d}\right)~=~
\left(R_{q}\right)^{ca}_{~db} \eqno (1.2) $$

As is well known, the quantum enveloping algebra is also a
Hopf algebra dual to ${\cal A}$, and the eq.(1.1) implies, in fact,
this dual relation (see Sec.II), so it must tightly connect with
the operators $L^{\pm}$. Ref.[11] gave their relations for
$U_{q}A_{1}$, but their generalization is not trivial. Some
authors tried to make clear the profound relations between
quantum groups and quantum enveloping algebras $^{[12]}$.
In this paper, after a further discussion on those relations
we carefully study the concrete expressions of the operators $L^{\pm}$
by the generators in a quantum enveloping algebra, and explicitly
give the expressions in two typical examples:
$U_{q}A_{N}$ and $U_{q}G_{2}$ algebras.

\vspace{10mm}
\section*{II. The dual relation}

\hskip 20pt
The generators $t_{j}$, $e_{j}$ and $f_{j}$ of a quantum enveloping
algebra $U_{q}{\cal G}$ satisfy the following algebraic relations:
$$\begin{array}{l}
t_{i}t_{j}~=~t_{j}t_{i},~~~~t_{j}^{-1}t_{j}~=~t_{j}t_{j}^{-1}~=~
{\bf 1},\\
t_{i}e_{j}~=~q_{i}^{a_{ij}}e_{j}t_{i},~~~~
t_{i}f_{j}~=~q_{i}^{-a_{ij}}f_{j}t_{i},~~~~q_{i}~=~q^{d_{i}}\\
{}~[e_{i},f_{j}]~=~\delta_{ij}\omega_{j}^{-1}\left(
t_{j}-t_{j}^{-1}\right),~~~~~~~~\omega_{j}~=~q_{i}-q_{i}^{-1}\\
\displaystyle \sum_{n=0}^{1-a_{ij}}~(-1)^{n}\left[\begin{array}{c}
1-a_{ij}\\[-2mm] n\end{array} \right]_{q_{i}}~e^{1-a_{ij}-n}
e_{j}e_{i}^{n}~=~0,~~~~i\neq j \\
\displaystyle \sum_{n=0}^{1-a_{ij}}~(-1)^{n}\left[\begin{array}{c}
1-a_{ij}\\[-2mm] n\end{array} \right]_{q_{i}}~f^{1-a_{ij}-n}
f_{j}f_{i}^{n}~=~0,~~~~i\neq j
\end{array} \eqno (2.1) $$

\noindent
where $a_{ij}$ is the Cartan matrix element of a Lie algebra
${\cal G}$, and $d_{j}$ is half-length of the simple root. For
the longer simple root, $d_{j}=1$. As a Hopf algebra, the coproduct
$\Delta$, the counit $\varepsilon$ and the antipode $S$ of the
generators are defined as follows:
$$\begin{array}{c}
\Delta(t_{j})~=~t_{j}\otimes t_{j},~~~~
\Delta(e_{j})~=~e_{j}\otimes {\bf 1}+t_{j}\otimes e_{j}, \\
\Delta(f_{j})~=~f_{j}\otimes t_{j}^{-1}+{\bf 1}\otimes f_{j}, \\
\varepsilon(t_{j})~=~1,~~~~\varepsilon(e_{j})~=~\varepsilon(f_{j})~=~0,\\
S(t_{j})~=~t_{j}^{-1},~~~~S(e_{j})~=~-t_{j}^{-1}e_{j},~~~~
S(f_{j})~=~-f_{j}t_{j},
\end{array} \eqno (2.2) $$

In fact, the representation matrix elements
$D^{[\lambda]}_{q}(\alpha)_{\mu \nu}$ of a quantum enveloping
algebra give the dual relation between the generators
$\alpha\in U_{q}{\cal G}$ and the non-commutative quantity
$\left(T^{[\lambda]}\right)^{\mu}_{~\nu}$ of the quantum group
in the corresponding representation:
$$\begin{array}{c}
\left\langle \alpha, \left(T^{[\lambda]}\right)^{\mu}_{~\nu}
\right\rangle~=~D^{[\lambda]}_{q}(\alpha)_{\mu \nu},~~~~
\alpha \in U_{q}{\cal G}
\end{array} \eqno (2.3) $$

\noindent
When $[\lambda]$ is the minimal representation,
$\left(T^{[\lambda]}\right)^{\mu}_{~\nu}$ becomes $T^{a}_{~b}$.

The universal ${\cal R}$ matrix can be expressed as
$${\cal R}~=~\displaystyle \sum_{i}~u^{i}\otimes v_{i}~\in~
U_{q}{\cal G}\otimes U_{q}{\cal G} \eqno (2.4) $$

\noindent
where $u^{i}$ and $v_{i}$ are the dual bases of up- and down-Borel
subalgebras, respectively. Take the value of ${\cal R}$ matrix in
the minimal representation, we get the $R_{q}$  matrix:
$$\langle {\cal R}, T^{c}_{~d}\otimes T^{a}_{~b}\rangle
{}~=~\left(R_{q}\right)^{ca}_{~db} \eqno (2.5) $$

The universal ${\cal R}$ matrix satisfies $^{[3,4]}$:
$${\cal R}\Delta(\alpha)~=~\left(P\circ \Delta(\alpha) \right){\cal R},
{}~~~~\alpha\in U_{q}{\cal G} \eqno (2.6) $$
$${\cal R}^{-1}~=~\left({\rm id}\otimes S^{-1}\right){\cal R}
\eqno (2.7) $$

\noindent
where  $P$ is the space transposition operator. Substituting eq.(2.6)
into eq.(2.5), we get
$$\begin{array}{l}
R^{ca}_{~rs}~\langle \Delta(\alpha), T^{r}_{~b}\otimes T^{s}_{~d}\rangle
{}~=~\langle {\cal R}\Delta(\alpha), T^{c}_{~b}\otimes T^{a}_{~d}\rangle \\
{}~=~\left\langle \left(P\circ \Delta(\alpha)\right){\cal R},
T^{c}_{~b}\otimes T^{a}_{~d}\right\rangle
{}~=~\langle \Delta(\alpha), T^{a}_{~s}\otimes T^{c}_{~r}\rangle
R^{rs}_{~bd} \end{array} \eqno (2.8) $$

\noindent
Due to arbitrariness of $\alpha$, one can immediately get eq.(1.1)
from eq.(2.8). Similarly, it is easy to show that equation (1.1)
holds for any other representation. So the consistent condition
(1.1) is just the direct reflection of the dual relation between
a quantum group and a quantum enveloping algebra.

Now we rewrite eq.(1.2) as the form of dual relation:
$$\left\langle \left(L^{-}\right)^{a}_{~b}, T^{c}_{~d}\right\rangle
{}~=~\left(R_{q}\right)^{ca}_{~db} $$

\noindent
Comparing it with eq.(2.5), we have
$$\left(L^{-}\right)^{a}_{~b}~=~
\langle {\cal R}, {\rm id} \otimes T^{a}_{~b}\rangle \eqno (2.9)$$

\noindent
Similarly,
$$\left(L^{+}\right)^{a}_{~b}~=~
\langle {\cal R}^{-1}, T^{a}_{~b}\otimes {\rm id}\rangle $$

\noindent
{}From eq.(2.7) we have:
$$\left(SL^{+}\right)^{a}_{~b}~=~
\langle {\cal R}, T^{a}_{~b}\otimes {\rm id}\rangle \eqno (2.10)$$

\noindent
In a different notation, the analogous relations were also given in
Ref.[13]. By making use of the general form of
the universal ${\cal R}$ matrix given by Jimbo $^{[3]}$, one can calculate
the explicit relations between $\left(L^{\pm}\right)^{a}_{~b}$ and
the generators in a quantum enveloping algebra from eqs.(2.9) and (2.10).

Denote by $Q_{+}$ the set of all the non-negative integral combinations
of simple roots,
$$Q_{+}~=~\left\{\beta|\beta=\displaystyle \sum_{j} m_{j} {\bf r}_{j},
m_{j}\in {\bar {\bf Z}_{-}} \right\} $$

\noindent
For a given $\beta \in Q_{+}$, we define two finite-dimensional
spaces $N_{\beta}^{\pm}$:
$$\begin{array}{ll}
N_{\beta}^{+}:~~~~&\left\{e_{i_{1}}e_{i_{2}}\cdots e_{i_{m}}|
\displaystyle \sum_{j}{\bf r}_{i_{j}}=\beta\right\}\\
N_{\beta}^{-}:~~~~&\left\{f_{i_{1}}f_{i_{2}}\cdots f_{i_{m}}|
\displaystyle \sum_{j}{\bf r}_{i_{j}}=\beta\right\}
\end{array} $$

\noindent
where ${\bf r}_{i_{j}}$ is the simple root corresponding to the
generators $e_{i_{j}}$ and $f_{i_{j}}$. Those two spaces both
are Hopf algebras, and they are associated and dual to each
other $^{[4]}$. In these dual spaces we can adopt the following dual
bases:
$$u^{i}\in N_{\beta}^{+},~~~~v_{i}\in N_{\beta}^{-},$$

\noindent
so that the general form of universal ${\cal R}$ matrix can be formally
 expressed as $^{[2,4]}$:
$${\cal R}~=~\left(\displaystyle \sum_{\beta \in Q_{+}}~L_{\beta} \right)
{}~q^{-H},~~~~
L_{\beta}~=~\displaystyle \sum_{i}~u^{i}\otimes v_{i}.
\eqno (2.11) $$

\noindent
A representation is called integrable if:

i) $V=\displaystyle \oplus_{\mu} V_{\mu}$
$$V_{\mu}~=~\left\{v_{\mu}\in V|t_{j}v_{\mu}=q^{\mu_{j}}_{j} v_{\mu},
{}~~\mu =\sum \mu_{j} \lambda_{j},~~1\leq j \leq N\right\}$$

ii) For any $v \in V$, there exists a common $M$ such that:
$$e_{j}^{M}~v~=~0,~~~~f_{j}^{M}~v~=~0. $$

\noindent
And $q^{-H}$ is given by
$$q^{-H} \left(v_{\mu}\otimes v_{\nu}\right)~=~
q^{-\displaystyle \sum_{ij}~d_{i}(a^{-1})_{ij}\mu_{i}\nu_{j}}
\left(v_{\mu}\otimes v_{\nu}\right) \eqno (2.12) $$

\noindent
where $(a^{-1})$ is the inverse Cartan matrix of the Lie
algebra ${\cal G}$.

When calculating $\left(L^{-}\right)^{a}_{~b}$ by eq.(2.9), the action
on the first subspace should keep its operator form, but for the second
subspace, it should be taken value in the minimal representation,
i.e., expressed as matrix form. By making use of this method, we
are able to express $\left(L^{\pm}\right)^{a}_{~b}$
operators by the generators in any quantum enveloping algebra in
principle. In the rest of this paper we will compute the explicit
expressions for two typical quantum enveloping algebras as examples.

\vspace{10mm}
\section*{III. $U_{q}A_{N}$ algebra}

\hskip 20pt
The reason for choosing $U_{q}A_{N}$ algebra as the first example
is that it is the simplest and the most useful qutantum enveloping
algebra. Jimbo $^{[2]}$ pointed out that for $U_{q}A_{N}$ any generator
corresponding to a non-simple root can be expressed as a
linear combination of the generators $e_{j}$ (or $f_{j}$) that
correspond to the simple roots:
$$\begin{array}{l}
E_{j(j+1)}~=~e_{j},~~~~E_{ij}~=~E_{ik} E_{kj}-q E_{kj} E_{ik},\\
E_{(j+1)j}~=~f_{j},~~~~E_{ji}~=~E_{jk} E_{ki}-q^{-1} E_{ki} E_{jk},
 \end{array},~~~~~~i<j<k \eqno (3.1) $$

\noindent
Through calculation, we find that $\left( L^{\pm} \right)^{a}_{~b}$ are
directly related with these generators.

In the minimal representation of $U_{q}A_{N}$, the bases can be
represented by only one index $a$, corresponding to the following
weight:
$$a \longrightarrow {\bf \lambda}_{a}-{\bf \lambda}_{a-1},~~~~
{\bf \lambda}_{0}={\bf \lambda}_{N+1}=0,~~~~1\leq a \leq N+1 \eqno (3.2) $$

\noindent
where ${\bf \lambda}_{a}$ are the fundamental dominant weights of
$U_{q}A_{N}$. The matrix forms of the generators in these bases are
given as follows:
$$\begin{array}{l}
t_{j}~|a\rangle~=~q^{\left(\delta_{aj}-\delta_{(a-1)j}\right)}~|a\rangle\\
e_{j}~|a\rangle~=~|(a+1)\rangle,~~~~
f_{j}~|a\rangle~=~|(a-1)\rangle  \end{array}~~~~1\leq a \leq N+1
 \eqno (3.3) $$

{}From eq.(2.12), we have
$$\begin{array}{l}
q^{-H}~\left(|a\rangle \otimes |b\rangle \right)~=~
\displaystyle \left\{\prod_{j}~t_{j}^{\tau_{jb}}
\otimes {\bf 1} \right\}~\left(|a\rangle \otimes |b\rangle \right) \\
\tau_{jb}~=~(a^{-1})_{j(b-1)}-(a^{-1})_{jb}~=~\left\{\begin{array}{ll}
\displaystyle {j \over N+1}-1,~~~~~&j\geq b \\
\displaystyle {j \over N+1},~~~~~&j< b \end{array} \right.
\end{array} $$

\noindent
or briefly express it as:
$$q^{-H}~\left(|a\rangle \otimes |b\rangle \right)~=~\left\{
\displaystyle \left(\prod_{j=1}^{N}~t_{j}^{j/(N+1)}\right)~
\displaystyle \left(\prod_{k=b}^{N}~t_{k}^{-1} \right)\otimes {\bf 1}
\right\}~\left(|a\rangle \otimes |b\rangle \right) \eqno (3.4) $$

Now we use eq.(2.9) to calculate $\left(L^{-}\right)^{a}_{~b}$. If $a<b$,
$\left(L^{-}\right)^{a}_{~b}$ vanishes. If
$a=b$, $\left(L^{-}\right)^{a}_{~b}$ only contains factor $q^{-H}$
given in eq.(3.4). If $a>b$, there is only one term in the
summation of ${\cal R}$ that has nonvanishing contribution to
$\left(L^{-}\right)^{a}_{~b}$, where $v_{i}$ is:
$$v_{i}~=~f_{a-1}~f_{a-2}~\cdots~f_{b},~~~~a\geq b \eqno (3.5) $$

\noindent
The key to the problem is to find $u^{i}$ dual to the above
$v_{i}$. In the calculation the following dual relations $^{[3,4]}$
are used:
$$\begin{array}{l}
\langle~t_{i}^{n}~,~t_{j}^{m}~\rangle~=~q_{i}^{-a_{ij}},~~~~
\langle~e_{i}~,~f_{j}~\rangle~=~- \delta_{ij}~\omega_{i}^{-1},\\
\langle~e_{i}~,~t_{j}^{m}~\rangle
{}~=~\langle~t_{i}^{n}~,~f_{j}~\rangle~=~0,
\end{array} \eqno (3.6) $$

\noindent
When $(a-b)$ is small, the dual operator $u^{i}$ is easy to
calculate. We can drew the general form of $u^{i}$ from these
$u^{i}$, then prove it by induction. Noting that for $U_{q}A_{N}$
all $\omega_{j}$ are equal to each other and can be denoted by
one symbol $\omega$. The result is:
$$\begin{array}{rl}
u^{i}&=~-\omega q^{a-b-1} \displaystyle \sum_{P}~
(-q)^{-n(P)}e_{p_{1}} e_{p_{2}} \cdots e_{p_{a-b}} \\
&=~(-1)^{a-b} \omega E_{ba} \end{array}~~~~~a>b \eqno (3.7) $$

\noindent
where $P$ is a certain permutation of $(a-b)$ objects, that can
be expressed by a product of some transpositions, and each
transposition moves a smaller number from the right of a bigger
one to the left.
$$P~=~\left(\begin{array}{cccc} (a-1)&(a-2)&\cdots&b\\
p_{1}&p_{2}&\cdots&p_{a-b} \end{array} \right) $$

\noindent
In such a multiplication expression, each transposition appears
at most once. We are only interested in the "neighboring
transposition" that permutes the places of two neighboring
numbers $c$ and $(c-1)$. Two $P$ are called equivalent if
their multiplication expressions contain same neighboring
transpositions without concerning their order.
The summation in eq.(3.6) runs over all inequivalent
permutations include the unit element. $n(P)$ is the
number of neighboring transpositions contained in $P$.

Substituting (3.4), (3.7) into (2.9) and (2.11), we finally get
$$\left(L^{-}\right)^{a}_{~b}~=~\left\{ \begin{array}{ll}
0,&a<b \\
\left\{\displaystyle \prod_{j=1}^{N} t_{j}^{j/(N+1)}\right\}
{}~\left\{\displaystyle \prod_{k=b}^{N} t_{k}^{-1}\right\}, &a=b \\
(-1)^{a-b} \omega E_{ba}
\left\{\displaystyle \prod_{j=1}^{N} t_{j}^{j/(N+1)}\right\}
{}~\left\{\displaystyle \prod_{k=b}^{N} t_{k}^{-1}\right\},~~~~&a>b
\end{array} \right. \eqno (3.8) $$

\noindent
Similarly
$$\left(L^{+}\right)^{a}_{~b}~=~\left\{ \begin{array}{ll}
(-1)^{b-a+1} \omega
\left\{\displaystyle \prod_{j=1}^{N} t_{j}^{-j/(N+1)}\right\}
{}~\left\{\displaystyle \prod_{k=a}^{N} t_{k}\right\} E_{ba},~~~~&a<b\\
\left\{\displaystyle \prod_{j=1}^{N} t_{j}^{-j/(N+1)}\right\}
{}~\left\{\displaystyle \prod_{k=a}^{N} t_{k}\right\}, &a=b \\
0, &a>b \end{array} \right. \eqno (3.9) $$

Note that it is only for $U_{q}A_{N}$ where the generators
corresponding to non-simple roots can be expressed as (3.1), so
that the expressions of $\left(L^{\pm}\right)^{a}_{~b}$
become so simple. For other algebras, the expressions become
much more complicated. $U_{q}G_{2}$ algebra is an example that
will be discussed in the next section.

\vspace{10mm}
\section*{IV. $U_{q}G_{2}$ algebra}

\hskip 20pt
$U_{q}G_{2}$ is another typical example where the lengths of
two simple roots are different. Let
$$\begin{array}{l}
q_{1}~=~q,~~~~q_{2}~=~q^{1/3},~~~~[m]~=~\displaystyle {q_{2}^{m}-q_{2}^{-m}
 \over q_{2}-q_{2}^{-1} },\\
\omega_{2}~=~q_{2}-q_{2}^{-1},~~~~\omega_{1}~=~q_{1}-q_{1}^{-1}
{}~=~[3]\omega_{2} \end{array}
\eqno (4.1) $$

\noindent
The minimal representation of $U_{q}G_{2}$ algebra is given as follows:
$$\begin{array}{l}
D_{q}(t_{1})~=~{\rm diag}\left\{1,~q,~q^{-1},1,~q,~q^{-1},~1\right\},\\
D_{q}(t_{2})~=~{\rm diag}\left\{q_{2},~q_{2}^{-1},~q_{2}^{2},1,~
q_{2}^{-2},~q_{2},~q_{2}^{-1}\right\},\\
D_{q}(e_{1})_{23}~=~D_{q}(f_{1})_{32}=
D_{q}(e_{1})_{56}~=~D_{q}(f_{1})_{65}=1,\\
D_{q}(e_{2})_{12}~=~D_{q}(f_{2})_{21}=
D_{q}(e_{2})_{67}~=~D_{q}(f_{2})_{76}=1,\\
D_{q}(e_{2})_{34}~=~D_{q}(f_{2})_{43}=
D_{q}(e_{2})_{45}~=~D_{q}(f_{2})_{54}=[2]^{1/2}
\end{array} \eqno (4.2) $$

\noindent
The rest of elements are vanishing. Through tedious calculation, we
obtain:
$$\begin{array}{ll}
\left(L^{-}\right)^{1}_{~1}~=~t_{1}^{-1}t_{2}^{-2},~~~~~
&\left(L^{-}\right)^{2}_{~2}~=~t_{1}^{-1}t_{2}^{-1},\\
\left(L^{-}\right)^{3}_{~3}~=~t_{2}^{-1},~~~~~
&\left(L^{-}\right)^{4}_{~4}~=~1,\\
\left(L^{-}\right)^{5}_{~5}~=~t_{2},~~~~~
&\left(L^{-}\right)^{6}_{~6}~=~t_{1}t_{2},\\
\left(L^{-}\right)^{7}_{~7}~=~t_{1}t_{2}^{2},~~~~~
&  \end{array} \eqno (4.3) $$
$$\begin{array}{ll}
\left(L^{-}\right)^{2}_{~1}~=~-\omega_{2} e_{2} t_{1}^{-1}t_{2}^{-2},~~~~~
&\left(L^{-}\right)^{3}_{~2}~=~-\omega_{1} e_{1} t_{1}^{-1}t_{2}^{-1},\\
\left(L^{-}\right)^{4}_{~3}~=~-\omega_{2} [2]^{1/2} e_{2} t_{2}^{-1},~~~~~
&\left(L^{-}\right)^{5}_{~4}~=~-\omega_{2} [2]^{1/2} e_{2} ,\\
\left(L^{-}\right)^{6}_{~5}~=~-\omega_{1} e_{1} t_{2},~~~~~
&\left(L^{-}\right)^{7}_{~6}~=~-\omega_{2} e_{2} t_{1}t_{2},
  \end{array} \eqno (4.4) $$
$$\begin{array}{l}
\left(L^{-}\right)^{3}_{~1}~=~\omega_{2}
\left(q e_{1} e_{2}-e_{2} e_{1}\right)t_{1}^{-1}t_{2}^{-2},\\
\left(L^{-}\right)^{4}_{~2}~=~- \omega_{2} [2]^{1/2}
\left(e_{1} e_{2}-q e_{2} e_{1}\right) t_{1}^{-1}t_{2}^{-1},\\
\left(L^{-}\right)^{5}_{~3}~=~ \omega_{2}^{2}
q_{2}^{-1} e_{2}^{2} t_{2}^{-1},\\
\left(L^{-}\right)^{6}_{~4}~=~\omega_{2} [2]^{1/2}
\left(q e_{1} e_{2}-e_{2} e_{1}\right) ,\\
\left(L^{-}\right)^{7}_{~5}~=~- \omega_{2}
\left(e_{1} e_{2}-q e_{2} e_{1}\right) t_{2},
  \end{array} \eqno (4.5) $$
$$\begin{array}{l}
\left(L^{-}\right)^{4}_{~1}~=~\omega_{2} q_{2}^{2} [2]^{-1/2}
\left\{e_{1}e_{2}^{2}-([6]/[3])e_{2}e_{1}e_{2}+e_{2}^{2}e_{1}\right\}
t_{1}^{-1}t_{2}^{-2},\\
\left(L^{-}\right)^{5}_{~2}~=~- \omega_{2} [2]^{-1}
\left\{e_{1}e_{2}^{2}-q_{2}^{2}[2]e_{2}e_{1}e_{2}
+q_{2}^{4}e_{2}^{2}e_{1}\right\}
t_{1}^{-1}t_{2}^{-1},\\
\left(L^{-}\right)^{6}_{~3}~=~- \omega_{2} [2]^{-1}
\left\{q_{2}^{4}e_{1}e_{2}^{2}-q_{2}^{2}[2]e_{2}e_{1}e_{2}
+e_{2}^{2}e_{1}\right\}
t_{2}^{-1},\\
\left(L^{-}\right)^{7}_{~4}~=~ \omega_{2} q_{2}^{2} [2]^{-1/2}
\left\{e_{1}e_{2}^{2}-([6]/[3])e_{2}e_{1}e_{2}+e_{2}^{2}e_{1}\right\},
  \end{array} \eqno (4.6) $$
$$\begin{array}{l}
\left(L^{-}\right)^{5}_{~1}~=~\omega_{2} q_{2} [2]^{-1}
\left\{e_{1}e_{2}^{3}
-([4]-q_{2}^{-1})e_{2}e_{1}e_{2}^{2}+q_{2}([4]-q_{2})
e_{2}^{2}e_{1}e_{2}-q_{2}e_{2}^{3}e_{1}\right\}
t_{1}^{-1}t_{2}^{-2},\\
\left(L^{-}\right)^{6}_{~2}~=~- \omega_{1}\omega_{2}q_{2}^{2}[6]^{-1}
\left\{e_{1}^{2}e_{2}^{2}+[3]e_{2}e_{1}^{2}e_{2}-([6][5]/[3][2])
e_{1}e_{2}^{2}e_{1}+e_{2}^{2}e_{1}^{2}\right)t_{1}^{-1}t_{2}^{-1},\\
\left(L^{-}\right)^{7}_{~3}~=~-\omega_{2} q_{2}[2]^{-1}
\left(q_{2}e_{1}e_{2}^{3}
-q_{2}([4]-q_{2})e_{2}e_{1}e_{2}^{2}+([4]-q_{2}^{-1})
e_{2}^{2}e_{1}e_{2}-e_{2}^{3}e_{1}\right\}t_{2}^{-1},
\end{array} \eqno (4.7) $$
$$\begin{array}{l}
\left(L^{-}\right)^{6}_{~1}~=~\omega_{1} q ([6][2])^{-1}
\left\{ q_{2}^{2}e_{1}^{2}e_{2}^{3}+q_{2}^{-2}e_{2}^{3}e_{1}^{2}
+q_{2}([4]-q_{2}^{-1})e_{2}e_{1}^{2}e_{2}^{2}\right. \\
{}~~~~~+q_{2}^{-1}([4]-q_{2})e_{2}^{2}e_{1}^{2}e_{2}
-q_{2}^{2}([4]-q_{2})([6]/[3])e_{1}e_{2}^{2}e_{1}e_{2}\\
{}~~~~~\left.-q_{2}^{-2}([4]-q_{2}^{-1})([6]/[3])e_{2}e_{1}e_{2}^{2}e_{1}
+([6][4]/[3][2])e_{1}e_{2}^{3}e_{1}
\right\}t_{1}^{-1}t_{2}^{-2},\\
\left(L^{-}\right)^{7}_{~2}~=~\omega_{1} q ([6][2])^{-1}
\left\{ q_{2}^{-2}e_{1}^{2}e_{2}^{3}+q_{2}^{2}e_{2}^{3}e_{1}^{2}
+q_{2}^{-1}([4]-q_{2})e_{2}e_{1}^{2}e_{2}^{2}\right. \\
{}~~~~~+q_{2}([4]-q_{2}^{-1})e_{2}^{2}e_{1}^{2}e_{2}
-q_{2}^{-2}([4]-q_{2}^{-1})([6]/[3])e_{1}e_{2}^{2}e_{1}e_{2}\\
{}~~~~~\left.-q_{2}^{2}([4]-q_{2})([6]/[3])e_{2}e_{1}e_{2}^{2}e_{1}
+([6][4]/[3][2])e_{1}e_{2}^{3}e_{1}
\right\}t_{1}^{-1}t_{2}^{-1},\\
\left(L^{-}\right)^{7}_{~1}~=~\omega_{1}\omega_{2}
q ([6][2])^{-1}\left\{([6][4]^{2}/[3][2])e_{2}e_{1}e_{2}^{2}e_{1}e_{2}
-([6]/[3][2])e_{1}e_{2}^{4}e_{1}\right.\\
{}~~~~~\left.-([4][2]-1)e_{2}^{2}e_{1}^{2}e_{2}^{2}
+(e_{2}^{4}e_{1}^{2}+e_{1}^{2}e_{2}^{4})
-([6][5]/[3][2])(e_{1}e_{2}^{2}e_{1}e_{2}^{2}
+e_{2}^{2}e_{1}e_{2}^{2}e_{1})
\right\}t_{1}^{-1}t_{2}^{-2},
 \end{array} \eqno (4.8) $$

$\left(L^{+}\right)^{a}_{~b}$ can be obtained from
$\left(L^{-}\right)^{b}_{~a}$  by following transformations: replace
$q_{j}$ by $q_{j}^{-1}$ (so $\omega_{j}$ becomes $-\omega_{j}$),
$e_{j}$ by $f_{j}$, $t_{i}$ by $t_{i}^{-1}$, and reverse the
product order of operators in each term. For example,
$$\left(L^{+}\right)^{1}_{~4}~=~-\omega_{2} q_{2}^{-2} [2]^{-1/2}
t_{1}t_{2}^{2}
\left\{f_{2}^{2}f_{1}-([6]/[3])f_{2}f_{1}f_{2}+f_{1}f_{2}^{2}\right\}$$

\noindent
It seems to us that this relation between $\left(L^{+}\right)^{a}_{~b}$
and  $\left(L^{-}\right)^{b}_{~a}$ holds for all quantum
enveloping algebras.

\vspace{10mm}
{\bf References}

\noindent
[1] V. G. Drinfel'd, Quantum Group, Proceedings of the International
Congress of Mathematicians, Berkeley, 1986, Vol.1, p.798.

\noindent
[2] M. Jimbo, Lett. Math. Phys. {\bf 10}(1985)63; Commun. Math. Phys.
{\bf 102}(1986)537.

\noindent
[3] M. Jimbo, Topics from representations of $U_{q}(g)$, - An
introductory guide to physicists, 1991, Nankai Lectures on
Mathematical Physics, World Scientific, Singapore, p.1.

\noindent
[4] Zhong-Qi Ma, Yang-Baxter Equation and Quantum Enveloping Algebras,
World Scientific, Singapore, 1993.

\noindent
[5] C. N. Yang, Phys. Rev. Lett. {\bf 19}(1967)1312; R. J.
Baxter, Ann. Phys. {\bf 70}(1972)193.

\noindent
[6] E. Date, M. Jimbo, K. Miki, and T. Miwa, Commun. Math.
Phys. {\bf 138}(1991)393.

\noindent
[7] M. Jimbo and T. Miwa, Algebraic Analysis of Solvable
Lattice Models, Conference Board of the Mathematical
Sciences, Regional Conference Series in Mathematics,
No. 85, 1995.

\noindent
[8] T. Brzezinski and S. Majid, Commun. Math. Phys. {\bf 157}(1993)591.

\noindent
[9] S. Watamura, Commun. Math. Phys. {\bf 158}(1993)67.

\noindent
[10] Bo-Yu Hou, Bo-Yuan Hou and Zhong-Qi Ma, J. Phys. A. {\bf 28}(1995)543.

\noindent
[11] L. D. Faddeev, N. Y. Reshetikhin and L. A. Takhtajan,
Quantization of Lie groups and Lie algebras, in Algebraic
Analysis, Academic Press, 1988, p.129.

\noindent
[12] Ke Wu, Han-Ying Guo and Ren-Jie Zhang, High Ener. Phys. and
Nucl. Phys. {\bf 17}(1993)262 (in chinese); Ke Wu and Han-Ying Guo, High
Ener. Phys. and Nucl. Phys. {\bf 17}(1993)796 (in chinese).

\noindent
[13] P. Schupp, P. Watts and B. Zumino, Commun. Math. Phys.
{\bf 157}(1993)305.

\end{document}